\documentclass{amsart}
\RequirePackage{url}
\usepackage{verbatim}
\usepackage {amsfonts}
\usepackage {   amsmath     }
\usepackage {   amssymb     }
\usepackage {   amsthm      }
\usepackage {   latexsym    }
\usepackage {   graphics    }
\usepackage {   color       }
\usepackage {   multicol    }

\usepackage {   amssymb     }
\usepackage{a4}

\theoremstyle{plain}

\theoremstyle{remark}

\newcommand{\mt}[4]{\left(\begin{array}{cc} #1 & #2 \\ #3 & #4\end{array}\right)}

\DeclareMathOperator{\ord}{ord}

\newcommand{\calO}{O}

\author{Mark van Hoeij$\mbox{}^{\dagger}$}\thanks{$\mbox{}^{\dagger}$ Supported by NSF 1618657.}
\author{Cristian-Silviu Radu$\mbox{}^{\ddagger}$}\thanks{$\mbox{}^{\ddagger}$ Supported by grant SFB F50-06 of the Austrian Science Fund (FWF)}

\begin{document}

\title{Computing an order complete basis for $M^{\infty}(N)$ and Applications}

\begin{abstract}
This paper gives a quick way to construct all modular functions for the group $\Gamma_0(N)$ having only a pole at $\tau = i \infty$.
We assume that we are given two modular functions $f,g$ for $\Gamma_0(N)$ with poles only at $i \infty$ and coprime pole orders.
As an application we obtain two new identities from which one can derive that $p(11n+6)\equiv 0\pmod{11}$, here $p(n)$ is the usual partition function. 
\end{abstract}
\maketitle

\section{Description of the Problem}
For basic notions about modular functions used in this paper we refer to \cite{PaRa2}.
In this paper we show how to obtain an order complete basis for $M^{\infty}(N)$ with an application to the case $N=11$.
We use this basis to obtain two new Ramanujan type identities for $\sum_{n=0}^{\infty}p(11n+6)q^n$.
Such bases have also been constructed by other authors \cite{Atkin,FGG,hughes,Kolberg,Kod,Lehner6,Ne1,Ne2}
by using various tricks to produce sufficiently many new modular functions $f_1,f_2,\ldots \in M^{\infty}(N)$ until
$\mathbb{C}[f_1,f_2,\ldots]$ becomes equal to $M^{\infty}(N)$.
The advantage of our approach is that we need only two functions in $t, f \in M^{\infty}(N)$.
Then $\mathbb{C}[t, f]$ will generally be a proper subset of $M^{\infty}(N)$, but instead of searching for more modular
functions, we fill this gap with a {\em normalized integral basis}.

Let $t$ and $f$ be modular functions for the group $\Gamma_0(N)$ with poles only at $\tau = i \infty$, in other words, let $t, f \in M^{\infty}(N)$.
Suppose that the pole orders are $n$ and $m$ respectively, and that $\gcd(n,m)=1$, such functions always exist \cite[Example 2.3]{PPRa}.
Then there exists an irreducible polynomial $p = p(x,y) \in\mathbb{C}[x,y]$ with $p(t, f)=0$, ${\rm deg}_x(p) = m$, and ${\rm deg}_y(p) = n$ by \cite[Lemma 1]{YY}.
One can compute $p$ from the $q$-expansions of $t$ and $f$ by making an Ansatz for the unknown coefficients of $p$
and solving a system of equations where each equation is a coefficient in the $q$-expansion of $p(t,f)$.
We use $p$ to compute in the function field $\mathbb{C}(t,f) \cong \mathbb{C}(x)[y]/(p)$.

The function field $\mathbb{C}(t,f)$ contains $M^{\infty}(N)$ see \cite[Prop 4.3]{PPRa}, here $M^{\infty}(N)$ is the set of all modular functions for the group $\Gamma_0(N)$ with a pole only at $i \infty$.
Obtaining all modular functions for the group $\Gamma_0(N)$ having a pole only at $i \infty$ is equivalent to finding all modular functions $h\in\mathbb{C}(t,f)$ that
are {\em integral} over $\mathbb{C}[t]$ (which means there is a {\em monic} polynomial $g(X) \in \mathbb{C}[t][X]$ for which $g(h)=0$).
Thus, one starts by computing an {\em integral basis}, which is a basis $b_1,\dots,b_n \in\mathbb{C}(t,f)$ of the $\mathbb{C}[t]$-module of all $h \in \mathbb{C}(t,f)$ that are integral over $\mathbb{C}[t]$.
There are several algorithms to compute an integral basis \cite{Montes, MvH} and implementations in several computer algebra systems.
Then every $h$ that is integral over $\mathbb{C}[t]$ can be written as $h=p_1(t)b_1+\dots+p_n(t)b_n$ for some polynomials $p_1,\ldots,p_n$.
Given the $q$-expansions of $h$ and $b_1,\ldots,b_n$ the algorithm described in \cite[Alg. MW]{Ra3} can find $p_1, \dots,p_n$
provided that $\ord_{i \infty}(b_1)<\ord_{i \infty}(b_2)<\dots<\ord_{i \infty}(b_n)$. We call such an integral basis {\em order complete}.

After computing an integral basis, we can find an order complete basis by using
{\em normalization at infinity} from Trager's PhD thesis \cite[Chapter 2, Section 3]{Trager}, see Section~\ref{Normalize} for details.

\subsection{Notations}
\mbox{} \\
$K = \mathbb{C}(x)[y]/(p)$ where $p \in \mathbb{C}[x,y]$ is irreducible. \\
${\calO}_K$ is ring of all elements of $K$ that are integral over $\mathbb{C}[x]$. \\
$R_{\infty}$ is the ring of all $h \in \mathbb{C}(x)$ that have no pole at $x=\infty$. \\
 ${\calO}_{\infty}$ is ring of all elements of $K$ that are integral over $R_{\infty}$.

To compute a basis of ${\calO}_{\infty}$ as $R_{\infty}$-module, first substitute
$x \mapsto 1/\tilde{x}$, then compute a {\em local integral basis} at $\tilde{x}=0$ (most
integral basis implementations allow the option of computing a local integral basis).
After that, replace $\tilde{x}$ by $1/x$.
  
\subsection{Normalize an integral basis at infinity}
\label{Normalize}

The process of normalizing an integral basis at infinity was introduced
in \cite{Trager} in order to compute a Riemann-Roch space that
was needed for integrating algebraic functions. For completeness we will
describe this process:

\noindent {\bf Algorithm:} Normalize an integral basis at infinity.
\begin{enumerate}
\item Let $b_1,\ldots,b_n$ be a basis of ${\calO}_K$ as $\mathbb{C}[x]$-module.
\item Let $b_1',\ldots,b_n'$ be a basis of ${\calO}_{\infty}$ as $R_{\infty}$-module.
\item Write $b_i = \sum_{j=1}^n r_{ij} b_i'$ with $r_{ij} \in \mathbb{C}(x)$.
\item Let $D \in \mathbb{C}[x]$ be a non-zero polynomial
for which $a_{ij} := D r_{ij} \in \mathbb{C}[x]$ for all $i,j$.
Now $D b_i = \sum_{j=1}^n a_{ij} b_i'$.
\item \label{goback}
For each $i \in \{1,\ldots,n\}$,
let $m_i$ be the maximum of the degrees of $a_{i1},\ldots,a_{in}$.
Now let $V_i \in \mathbb{C}^n$ be the vector whose $j$'th entry is
the $x^{m_i}$-coefficient of $a_{ij}$.
Let $d_i := m_i - {\rm deg}_x(D)$.
\item If $V_1,\ldots,V_n$ are linearly independent, then
return $b_1,\ldots,b_n$ and $d_1,\ldots,d_n$ and stop. \\
Otherwise, take $c_1,\ldots,c_n \in \mathbb{C}$, not all 0, for which
$c_1 V_1 + \cdots c_n V_n = 0$.
\item Among those $i \in \{1,\ldots,n\}$ for which $c_i \neq 0$,
choose one for which $d_i$ is maximal. For this $i$, do the
following
\begin{enumerate}
\item \label{newb} Replace $b_i$ by $\sum_{k=1}^n c_k x^{d_i - d_k} b_k$.
\item Replace $a_{ij}$ by $\sum_{k=1}^n c_k x^{d_i - d_k} a_{kj}$ for
all $j \in \{1,\ldots,n\}$.
\end{enumerate}
\item Go back to step~\ref{goback}.
\end{enumerate}
The $b_1,\ldots,b_n$ remain a basis of ${\calO}_K$ throughout
the algorithm because the new $b_i$ in step~\ref{newb} can be written
as a nonzero constant times the old $b_i$ plus
a $\mathbb{C}[x]$-linear combination of the $b_j$, $j\neq i$.
When we go back to step~\ref{goback} the non-negative integer $d_i$
decreases while the $d_j$, $j \neq i$ stay the same. Hence the
algorithm must terminate.

Let $b_1,\ldots,b_n$ and $d_1,\ldots,d_n$ be the output of the algorithm.
By construction, the number $d_i$ in the algorithm is the smallest
integer for which $b_i \in x^{d_i} {\calO}_{\infty}$.
If $\beta \in {\calO}_K$ with $\beta \neq 0$ then we can write
$\beta = c_1 b_1 + \cdots + c_n b_n$ for some
$c_1,\ldots,c_n \in \mathbb{C}[x]$.
Denote $d_{\beta}$ as the maximum of ${\rm deg}_x(c_j)+d_j$ taken
over all $j$ for which $c_j \neq 0$.
Then $\beta \in x^{d_{\beta}} {\calO}_{\infty}$
by construction. Since the vectors $V_1,\ldots,V_n$ in
the algorithm are linearly independent when the algorithm terminates,
there can not be any cancellation, which means that ${d_{\beta}}$
is the smallest integer for which
$\beta \in x^{d_{\beta}} {\calO}_{\infty}$.
Because of this, we get the following:

\noindent  If $d$ is a positive integer, then the set
$B_d := \{ x^j b_i \ | \ 0 \leq j \leq d-d_i, \ 1\leq i \leq n\}$ is a basis of
${\calO}_K \bigcap x^d {\calO}_{\infty}$ as $\mathbb{C}$-vector space.

Note that $B_d$ is a basis of the Riemann-Roch
space of the pole-divisor of $x^d$. So computing $B_d$ can be interpreted
as (i): a direct application of a normalized integral basis, or (ii): a
special case of algorithms \cite{RiemannRochSpaces,RiemannRochSpaces2} for Riemann-Roch spaces.
The two interpretations are equivalent because the first step in
computing Riemann-Roch spaces is to compute a normalized
integral basis.

We can take $q$-expansions for each of the elements of $B_d$, and then make a
change of basis so that the new basis $B_d^{\rm REF}$ will have $q$-expansions in Reduced Echelon Form.
This means that if $b \in B_d^{\rm REF}$ and $b = a_r q^r + a_{r+1} q^{r+1} + \cdots$ with $a_r \neq 0$ then
$a_r = 1$ and all other basis elements have a zero coefficient at $q^r$.
Then $B_d^{\rm REF}$, for suitable $d$, is an order complete basis. For an implementation and two
examples see:   \\
www.math.fsu.edu/$\sim$hoeij/files/OrderComplete

\section{New Identities}
We will give two identities of Ramanujan type found using our algorithm (the second one is only on our website).
Let $p(n)$ be the partition function. Define
$$t:=q^{-5} \prod_{n=1}^{\infty}\left(\frac{1-q^n}{1-q^{11n}}\right)^{12}.$$
and
$$h:=qt\prod_{n=1}^{\infty}(1-q^{11k}) \sum_{n=0}^{\infty}p(11n+6)q^n$$
and
$$f:=(dt/dq)\prod_{n=1}^{\infty}(1-q^n)^{-2}(1-q^{11n})^{-2}.$$
Both $h$ and $t$ are modular functions in $M^{\infty}(11)$, see \cite[Lemma 3.1]{PaRa2}.

To prove that $f$ is in $M^{\infty}(11)$ as well,
first note that by \cite[Prop. 3.1.1]{Ligozat}
$$b(\tau):=q\prod_{n=1}^{\infty}(1-q^n)^2(1-q^{11n})^2,\quad q=e^{2\pi i\tau}$$
satisfies
\begin{equation}\label{rell1}
  b\Bigl(\frac{a\tau+b}{c\tau+d}\Bigr)=(c\tau+d)^{2}b(\tau)
  \end{equation}
for all $\mt{a}{b}{c}{d}\in \Gamma_0(11)$. 
Since $t\in M^{\infty}(11)$, we have
$$t\Bigl(\frac{a\tau+b}{c\tau+d}\Bigr)=t(\tau).$$
The derivative with respect to $\tau$ is:
\begin{equation}\label{rell2}
  (c\tau+d)^{-2}t'\Bigl(\frac{a\tau+b}{c\tau+d}\Bigr)=t'(\tau)
\end{equation}
Multiplying  \eqref{rell2} by $(c\tau+d)^2$ and dividing by \eqref{rell1} gives
$$(t'/b)\Bigl(\frac{a\tau+b}{c\tau+d}\Bigr)=(t'/b)(\tau).$$
Since $\frac{d}{d\tau}=2\pi iq\frac{d}{dq}$, it follows that $t'/b=2\pi if$.  
Therefore
$$f\Bigl(\frac{a\tau+b}{c\tau+d}\Bigr)=f(\tau)$$
for all $\mt{a}{b}{c}{d}\in \Gamma_0(11)$. 
Furthermore, since $b(\tau)$ has no zeros in the upper half plane and $t(\tau)$ is holomorphic in the upper half plane it follows that $f$ is holomorphic in the upper half plane. Hence the first condition of being a modular function for $\Gamma_0(11)$ according to the definition in \cite{PaRa2} is satisfied. The second condition is equivalent to showing that for $\gamma=\mt{a}{b}{c}{d}\in \mathrm{SL}_2(\mathbb{Z})$ we have an expansion of the form
\begin{equation}\label{thisprop}
  f\Bigl(\frac{a\tau+b}{c\tau+d}\Bigr)=\sum_{n=m(\gamma)}^{\infty}a_{\gamma}(n)q^{\frac{\gcd(c^2,n)n}{N}}.
  \end{equation}
As seen in \cite{PaRa2}, if this property hold for $\mt{a}{b}{c}{d}$, then it also holds for $\mt{a'}{b'}{c'}{d'}$, if there exists $\mt{A}{B}{C}{D}\in \Gamma_0(11)$ such that $\frac{A\frac{a}{c}+B}{C\frac{a}{c}+D}=\frac{a'}{c'}$.
So we need to find representatives of the orbits of the action of $\Gamma_0(11)$ on $\mathbb{Q}\cup \{i \infty\}$, that is, the cusps of $\Gamma_0(N)$. From \cite{RaEq} we find that these representatives are $0$ and $i \infty$. 
Then it suffices to show \eqref{thisprop} for two cases: $\mt{1}{0}{0}{1}$ and $\mt{0}{-1}{1}{0}$. The first case holds because $f$ is a $q$-series. For the second case we need to show that $f(-1/\tau)$ is a Laurent series in $q^{1/11}$ with finite principal part. By \cite{Ra} we have
$$\eta(-1/\tau)=(-i\tau)^{1/2}\eta(\tau).$$
This implies
\begin{equation}\label{rrid}
  t(-1/\tau)=t^{-1}\Bigl(\frac{\tau}{11}\Bigr)
  \end{equation}
and
$$b(-1/\tau)=-\frac{1}{11}b(\frac{\tau}{11})\tau^2.$$
The derivative of \eqref{rrid} is
$$\tau^{-2}t'(-1/\tau)=\frac{1}{11}t^{-2}(\frac{\tau}{11})t'(\frac{\tau}{11})$$
which is equivalent to
$$t'(-1/\tau)=\frac{1}{11}\tau^2t^{-2}(\frac{\tau}{11})t'(\frac{\tau}{11}).$$
This implies
$$(t'/b)(-1/\tau)=-(t'/b)(\frac{\tau}{11})t^{-2}(\frac{\tau}{11}).$$
Hence
$$f(-1/\tau)=-f(\tau/11)t^{-2}(\tau/11)=5q^{4/11}+O(q^{5/11}).$$
So the last condition for $f$ being a modular function for $\Gamma_0(11)$ is verified. In order for $f$ to be in $M^{\infty}(11)$ we need the order of $f$ to be nonnegative at all cusps except $i \infty$.
That only leaves the cusp $0$ where the order is $4$. This shows $f \in M^{\infty}(11)$.

We want to express $h$ as an element of $\mathbb{C}(t,f)$.
The pole orders of $t$ and $f$ are $5$ and $6$ so $p(x,y) = \sum_{i=0}^6 \sum_{j=0}^5 a_{ij} x^i y^j$ is an Ansatz for the algebraic relation $p(t,f)=0$.
Solving linear equations coming from $q$-expansions gives $$p(x,y) = 
y^5+170xy^4+9345x^2y^3+167320x^3y^2+(5^5x^2-7903458x+5^5 11^6)x^4.$$
We use $p(x,y)$ to compute in $\mathbb{C}(t,f) \cong \mathbb{C}(x)[y]/(p)$.
We compute $B^{\rm REF}_d$ from the previous section with $d=1$
and obtain $b_0, b_2, b_3, b_4, b_5$ where $b_0=1$ and $b_i = q^{-i} + c_i q^{-1} + O(q^1)$ for $i = 2,\ldots,5$ for some constants $c_i$.
Since $h$ has a pole of order 4, we can write it as a linear combination of $b_0, b_2, b_3, b_4$.
We have $b_0 = 1$ and
\begin{eqnarray*}
b_2   \, =& \hspace{-18pt} 12   +   \frac{5t}{22} \left(  \ \ \ \frac{t-11^3}{f+47t}  \, - \,  \frac{(42t+  f)(t+11^3)}{f^2+89ft+1424t^2}\right)   & \hspace{-6pt} = \ q^{-2}+2q^{-1}+5q+8q^2+O(q^3) \\
b_3   \, =& \hspace{-12pt} 12   +  \frac{5t}{22}\left( \ \, 3 \frac{t-11^3}{f+47t} -   \frac{(16t+3f)(t+11^3)}{f^2+89ft+1424t^2}\right)   & \hspace{-6pt} = \ q^{-3}+ \ q^{-1} \, +2q+2q^2+O(q^3) \\
b_4   \, =& \hspace{-4pt} 12   +   \frac{5t}{22}\left(-3 \frac{t-11^3}{f+47t} - \frac{(28 t + 19 f)(t+11^3)}{f^2+89ft+1424t^2}\right)  & \hspace{-6pt} = \ q^{-4}-2q^{-1}+6q+3q^2+O(q^3). \\
\end{eqnarray*}
Like in \cite[Alg. MW]{Ra3}, we use
$$h=11q^{-4}+165q^{-3}+748q^{-2}+1639q^{-1}+3553+O(q)$$
to find
$$h-11b_4-165b_3-748b_2-3553b_0=O(q).$$
This expression in $M^{\infty}(11)$ has no poles and a root at $\tau = i \infty$ (at $q=0$) hence it is the zero function. Therefore
$$h=11b_4+165b_3+748b_2+3553b_0.$$

Replacing $b_0,b_2,b_3,b_4$ with their corresponding expressions in terms of $t$ and $f$ gives
$$h    = qt\prod_{n=1}^{\infty}(1-q^{11k}) \sum_{n=0}^{\infty}p(11n+6)q^n
       = 11^4 + 55t \left( 5\frac{t-11^3}{f+47t} - \frac{2(71t+3f)(t+11^3)}{f^2+89ft+1424t^2} \right).$$


This implies $p(11n+6)\equiv 0\pmod{11}$. Other expressions for $h$ that prove this congruence
were already in \cite{Atkin, Lehner6}, however, our expression in terms of $t,f$ is novel.

For our second example, take $t$ and $h$ be as before and let
$$E_4:=1+240\sum_{n=1}^{\infty}\frac{n^3q^n}{1-q^n}$$
be the usual Eisenstein series. Let
$$\Delta:=q\prod_{n=1}^{\infty}(1-q^n)^{24}.$$
Let $J:=E_4^3/\Delta= q^{-1}+\cdots$ and
$$f:=Jt^3.$$

Next we show that $f\in M^{\infty}(11)$. From the last chapter of \cite{Se} we find
$$E_4\Bigl(\frac{a\tau+b}{c\tau+d}\Bigr)=(c\tau+d)^4E_4(\tau)$$
and
$$\Delta\Bigl(\frac{a\tau+b}{c\tau+d}\Bigr)=(c\tau+d)^{12}\Delta(\tau)$$
for all $\mt{a}{b}{c}{d}\in\mathrm{SL}_2(\mathbb{Z})$. These two identities imply
$$J\Bigl(\frac{a\tau+b}{c\tau+d}\Bigr)=J(\tau)$$
for all $\mt{a}{b}{c}{d}\in \mathrm{SL}_2(\mathbb{Z})$. 
Since $\mathrm{SL}_2(\mathbb{Z})$ has only one cusp, $i \infty$, and since $J$ is a $q$-series it follows that $J$ is a modular function on $\mathrm{SL}_2(\mathbb{Z})$ and thus on $\Gamma_0(11)$.

Since $t(\tau)$ is already a modular function on $\Gamma_0(11)$, it follows that $f$ is a modular function on $\Gamma_0(11)$. To show that $f$ is in $M^{\infty}(11)$ it suffices to show that the order of $f$ at the cusp $0$ is nonnegative.
Since $J(-1/\tau)=(q^{-1/11})^{11}+O(1)$ the order of $J$ at $0$ is $-11$. The order of $t$ at $0$ is $5$, so the order of $f$ at the cusp $0$ is $-11+3 \cdot 5 = 4 \geq 0$.  This shows $f\in M^{\infty}(11)$. 

The only pole of $f$ is at $i \infty$, it has order 16. We compute the algebraic relation $p(t,f)=0$ with the Ansatz method, and use $p$ to compute $B_d^{\rm REF}$.
Then we express $h$ in terms of the $t$ and the new $f$. This relation, and the Maple file that computes it, are given at www.math.fsu.edu/$\sim$hoeij/files/OrderComplete.


\bibliography{zzz2}
\bibliographystyle{plain}

\end{document}